\documentclass[11pt]{amsart}
\usepackage{amssymb,amsmath,txfonts,mathrsfs}
\newtheorem{theorem}{Theorem}[section]

\newtheorem{lemma}[theorem]{Lemma}
\newtheorem{remark}[theorem]{Remark}
\newtheorem{claim}[theorem]{Claim}
\newtheorem{note}[theorem]{Note}

\newtheorem{definition}[theorem]{Definition}
\newtheorem{cor}[theorem]{Corollary}

\numberwithin{equation}{section}

\def\pf{{\it Proof:}~}

\begin{document}

\title[The short time asymptotics of Nash entropy]{The short time asymptotics of Nash entropy}
\author{Guoyi Xu}
\address{Guoyi Xu\\ Mathematics Department\\University of California, Irvine\\CA 92697\\United States}
\email{guoyixu@math.uci.edu}
\date{\today}

\begin{abstract}

Let $(M^n, g)$ be a complete Riemannian manifold with $Rc\geq -Kg$, $H(x, y, t)$ is the heat kernel on $M^n$, and $H= (4\pi t)^{-\frac{n}{2}}e^{-f}$. Nash entropy is defined as $N(H, t)= \int_{M^n} (fH) d\mu(x)- \frac{n}{2}$.  We studied the asymptotic behavior of $N(H, t)$ and $\frac{\partial}{\partial t}\Big[N(H, t)\Big]$ as $t\rightarrow 0^{+}$, and got the asymptotic formulas at $t= 0$. In the Appendix, we got Hamilton-type upper bound for the Laplacian of positive solution to the heat equation on such manifolds, which has its own interest.
\\[3mm]
\end{abstract}

\maketitle

\section{Introduction}
On a complete manifold $(M^n, g)$ with $Rc\geq -Kg$, where $K> 0$ is a consant, for fixed $y\in M^n$, it is well-known that the heat kernel $H(x, y, t)$ on $(M^n, g)$ is unique. We assume $H= (4\pi t)^{-\frac{n}{2}}e^{-f}$. As in \cite{Niadd}, Nash entropy is defined as: 
\begin{definition}\label{def Nash}
{\begin{equation}\label{Nash}
N(H, t)= \int_{M^n} (fH) d\mu(x)- \frac{n}{2}
\end{equation}
}
\end{definition}
Nash entropy has close relation to $\mathscr{W}$-entropy for linear heat equation, and the large time asymptotics of this entropy reflects the volume growth rate of the manifold (see \cite{Ni}, \cite{Niadd} and \cite{Nilarge}). 

In this paper, we studied the asymptotic behavior of $N(H, t)$ and $\frac{\partial}{\partial t}N(H, t)$ as $t\rightarrow 0^{+}$, and solved one problem proposed in \cite{RFTA} (Problem $23.36$ there). More precisely, we proved the following theorem:
\begin{theorem}\label{thm 3.2}
{Let $(M^n, g)$ be a complete Riemannian manifold with $Rc\geq -Kg$, where $K> 0$ is a constant. Then 
\begin{equation}\label{3.6}
{N(H, t)= -\frac{1}{2}R(y)\cdot t+ O(t^{\frac{3}{2}})
}
\end{equation}
and 
\begin{equation}\label{3.6.1}
{\frac{\partial}{\partial t}\Big[N(H, t) \Big]= -\frac{1}{2}R(y)+ o(1)
}
\end{equation}
where $\limsup_{t\rightarrow 0}O(t^{\frac{3}{2}})t^{-\frac{3}{2}}$ is bounded, $\lim_{t\rightarrow 0}o(1)= 0$ and $t$ is small enough.
}
\end{theorem}

One motivation to study the short time asymptotics of Nash entropy is Li-Yau-Perelman type estimate for heat equation on manifolds with Ricci curvature bounded from below. Motivated by Perelman's differential Harnack estimate for Ricci flow, in \cite{Ni}, on a closed manifold $(M^n, g)$ with $Rc\geq 0$, Ni proved the following Li-Yau-Perelman type estimate for the heat equation when $t> 0$:
\begin{equation}\label{0.2.0}
{2\Delta f(x, y, t)- |\nabla f(x, y, t)|^2+ \frac{f(x, y, t)- n}{t}\leq 0
}
\end{equation}
where $H(x, y, t)= (4\pi t)^{-\frac{n}{2}}e^{-f}$ is the heat kernel. In fact, (\ref{0.2.0}) is also true for heat kernel on complete manifold $(M^n, g)$ with $Rc\geq 0$ (see \cite{RFAA}). 

In the well-known paper \cite{Pere}, Perelman made the following claim (see remark $9.6$ there):
\begin{claim}\label{claim}
{If $(M^n, g)$ is a compact Riemannian manifold, $g_{ij}(x, t)$ evolves according to $\Big(g_{ij}\Big)_t= A_{ij}(t)$ and $g_{ij}(x, 0)= g_{ij}(x)$, $t\in (-T, 0]$. Define $\square = \frac{\partial}{\partial t}- \Delta$ and its conjugate $\square^{*}= -\frac{\partial}{\partial t}- \Delta- \frac{1}{2}A$ (where $A= g^{ij}A_{ij}$). Consider the fundamental solution $u= (-4\pi t)^{-\frac{n}{2}} e^{-f}$ for $\square^{*}$, starting as $\delta$-function at some point $(p, 0)$. Then for general $A_{ij}$ the function $\Big(\square \bar{f}+ \frac{\bar{f}}{t}\Big)(q, t)$, where $\bar{f}= f- \int_{M^n} fu$, is of order $O(1)$ for $(q, t)$ near $(p, 0)$.
}
\end{claim}

We will focus on the special case where the evolving metrics are the static metric. From Theorem \ref{thm 3.2}, it is easy to show that Perelman's claim in the static metric case is equivalent to the following claim on compact manifolds:
\begin{equation}\label{0.2}
{2\Delta f(x, y, t)- |\nabla f(x, y, t)|^2+ \frac{f(x, y, t)- n}{t}= -R(y)+ O(t+ d^2(x, y))
}
\end{equation}
If (\ref{0.2}) is true, it will be an improvement of (\ref{0.2.0}) when $t+ d^2(x, y)$ is small enough and $R(y)> 0$. But using the following explicit formula of heat kernel on hyperbolic manifold $\mathbb{H}^3$ (cf. section $9.2$ in \cite{Gri}):
\begin{equation}\nonumber
{H= (4\pi t)^{-\frac{3}{2}}\frac{d}{\sinh d}\exp{\big(-\frac{d^2}{4t}- t\big)}
}
\end{equation}
it is easy to check that (\ref{0.2}) is not true generally. Hence Claim \ref{claim} is not generally true for static metric case on complete manifolds. 

As observed in \cite{Niadd}, the integrand of $\frac{\partial}{\partial t}\Big[N(H, t)\Big]$ is nothing but the expression in Li-Yau's gradient estimate for heat kernel multiplying with the heat kernel, which is $-(\Delta \ln H+ \frac{n}{2t})H$. Because so far there is no sharp Li-Yau-type gradient estimate for heat kernel or solutions to the heat equation on complete manifolds with Ricci curvature bounded from below by negative constant, we hope that (\ref{3.6.1}) will be helpful on understanding this estimate better.

On the other hand, when $(M^n, g)$ is a compact Riemannian manifold, the short time behavior of logarithm of heat kernel had been studied by many probabilists. Although the heat kernel $H(x, y, t)$ has infinite sequence expansion at $t= 0$, generally there is no such expansion of $\ln H$ at $t= 0$, and the singularity of $\ln H$ at $t= 0$ can have many complicated situations. However, in \cite{Va}, Varadhan proved 
\begin{equation}\label{0.1}
{\lim_{t\rightarrow 0} t\ln H(x, y, t)= -\frac{d^2(x, y)}{4}
}
\end{equation}

Moreover, using stochastic processes methods, Malliavin and Stroock proved that the above equation is preserved while taking the first and second spatial derivatives on domain outside of cut locus (see \cite{MS}). Using analytic methods, (\ref{0.1}) is proved for complete Riemannain manifolds in \cite{CLY} by Cheng, Li and Yau. We hope that Theorem \ref{thm 3.2} will be useful on studying the short time behavior of logarithm of the heat kernel on complete manifolds by analytic methods.

The strategy to prove (\ref{3.6}) is using infinite sequence expansion $H_{N}(x, y, t)$ of $H(x ,y ,t)$ at $t= 0$, although generally $\ln H_{N}$ does not converge to $\ln H$ near $t= 0$ uniformly. In integral sense of (\ref{Nash}), we show there is a uniform convergence in Lemma \ref{lem 3.1} by using an improved estimate of $H- H_{N}$ got in Theorem \ref{thm 2.2}. The rest calculation about integral of $H_{N}$ is standard, for completeness we give details in full. 

To prove (\ref{3.6.1}), because the manifold $M^n$ can be non-compact, we need to be more careful on the switch of the order of differentiation and integration. The detailed proof of the validity of the switch is given in the beginning of section $4$. We need an upper bound of $\frac{H_{t}}{H}$ in verifying the above switch. This type bound is known for closed manifolds from \cite{Ham}, and in \cite{RFAA} (also see \cite{LYH}) the proof is sketched for complete manifolds with $Rc\geq 0$ following similar strategy of Kotschwar in \cite{Kots}. The detailed proof of this Hamilton-type upper bound for complete manifolds with $Rc\geq -Kg$ is included in the Appendix for completeness. 

\begin{note}\label{note 1.1}
After the paper is circulated and posted on the arXiv, Jia-yong Wu kindly informed us that he had independently proved Hamilton-type upper bound for complete manifolds with $Rc\geq 0$ in details, which is in \cite{Wu}. 
\end{note}


The paper is organized as follows: In section $2$, we state some preliminary results about heat kernel and get some improved estimates about $H-H_{N}$. In section $3$, we prove (\ref{3.6}). Using (\ref{3.6}) and results in Appendix, (\ref{3.6.1}) is proved in section $4$. In Appendix, on complete manifolds with Ricci curvature bounded from below, Hamilton-type upper bound of $\frac{H_{t}}{H}$ is proved. 

Acknowledgement: The author would like to thank Zhiqin Lu, Brett Kotschwar for interest and suggestions, and Peter Li, Jiaping Wang for their interest.

\section{Preliminary}
We firstly define some notations and functions. In the rest of the paper, we fix $y\in M^n$, define \[\Omega_y=\{x\in M^n: d(x, y)< inj_{g}(y)\}\]
where $inj_{g}(y)$ denotes the injectivity radius of metric $g$ at $y$. Define 
\[B(\rho)= \{x|\ d(x, y)\leq \rho\} \quad and \quad B_{z}(\rho)= \{x| \ d(x, z)\leq \rho \}\] 
Hence $B(\rho)= B_y(\rho)$. $V(B_z(\rho))$ is used to denote the volume of $B_z(\rho)$, $V_{-K}(\rho)$ is the volume of the geodesic ball of radius $\rho$ in the constant $\Big(-\frac{K}{n -1}\Big)$ sectional curvature space form. 

Choose $r\in (0, \frac{1}{4}inj_{g}(y))$, fix it and let $N_{0}= \frac{n}{2}+ 3$. Define $E= (4\pi t)^{-\frac{n}{2}}\exp \Big(-\frac{d^2(x, y)}{4t}\Big)$ and $\tilde{E}= (4\pi t)^{-\frac{n}{2}}\exp \Big(-\frac{d^2(x, y)}{5t}\Big)$. Sometimes we will use $B$ as the simplification of the notation $B(\frac{r}{2})$, $d(x, y)$ will be simplified as $d$, and what it means will be clear from the context. 

Assume $\eta: [0, \infty)\rightarrow [0, 1]$ is a $C^{\infty}$ cut-off function with 
\begin{equation}\label{2.1}
\eta(s)=\left\{
\begin{array}{rl}
&1\quad if\ s\leq r \\
&0\quad if\ s\geq 2r 
\end{array} \right.
\end{equation}
The following theorem collects some known results about heat kernel on complete manifolds (see \cite{RFTA}, \cite{GL}, \cite{Li} etc.).

\begin{theorem}\label{thm 2.1}
{$(M^n, g)$ is a complete Riemannian manifold with $Rc\geq -Kg$, where $K> 0$ is a constant. Then there exists a unique positive fundamental solution $H(x, y, t)$ to the heat equation, which is called the heat kernel. Moreover $H(x, y, t)\in C^{\infty}(M^n\times M^n\times (0, \infty))$ is symmetric in $x$ and $y$, and

$(i)$ \begin{equation}\label{2.3}
{\int_{M^n} H(x, y, t)d\mu(x)\equiv 1
}
\end{equation}
$(ii)$ \begin{equation}\label{2.4}
{H(x, y, t)= P_{N_0}(x, y, t)+ F_{N_0}(x, y, t)
}
\end{equation}

\begin{equation}\label{2.5}
{P_{N_0}(x, y, t)= \eta(d(x, y)) H_{N_0}(x, y, t)
}
\end{equation}
and 
\begin{equation}\label{2.6}
{H_{N_0}(x, y, t)= (4\pi t)^{-\frac{n}{2}}\exp \Big(-\frac{d^2(x, y)}{4t}\Big)\cdot \sum_{k= 0}^{N_0} \varphi_{k}(x, y)t^k
}
\end{equation}
$\varphi_{k}(x,y )\in C^{\infty}(\Omega_{y})$, $k= 0, 1, \cdots, N_0$. Also $H_{N_0}$ satisfies the following:
\begin{equation}\label{2.7}
{(\Delta- \frac{\partial}{\partial t})H_{N_0}(x, y, t)= E\Delta \varphi_{N_0} t^{N_0}
}
\end{equation}
$(iii)$ Let $\{x^k\}_{k=1}^n$ be exponential normal coordinates centered at $y\in M^n$, then $\varphi_{0}$ and $\varphi_{1}$ have the following asymptotic expansion:
\begin{equation}\label{2.9}
{\varphi_{0}(x, y)= 1+ \frac{1}{12}R_{pq}(y)x^px^q+ O(d^3(x, y))
}
\end{equation}
\begin{equation}\label{2.10}
{\varphi_{1}(x, y)= \frac{R(y)}{6}+ O(d(x, y))
}
\end{equation}
}
\end{theorem}

We will prove an estimate for $F_{N_0}$, this estimate is an improvement of the usual estimate of $F_{N_0}$, which only gives $t^{N_0+ 1- \frac{n}{2}}$ bound. The improved estimate (\ref{2.11}) is the key to the proof of Lemma \ref{lem 3.1}.

\begin{theorem}\label{thm 2.2}
{For $F_{N_0}(x, y, t)$ in Theorem \ref{thm 2.1}, we have the following estimates:
\begin{equation}\label{2.11}
{|F_{N_0}(x, y, t)|\leq Ct^4\exp{(-\frac{d^2(x, y)}{5t})}
}
\end{equation}
and
\begin{equation}\label{2.11.0}
{\Big|\frac{\partial}{\partial t}F_{N_0}(x, y, t)\Big|\leq Ct^2\exp{(-\frac{d^2(x, y)}{5t})}
}
\end{equation}
where $t$ is small enough and $C$ is a positive constant independent of $x$, $t$. 
}
\end{theorem}

\begin{remark}\label{remark 2.1}
{(\ref{2.11}) had been proved in \cite{GL} for uniformly parabolic operators. Our proof of (\ref{2.11}) and (\ref{2.11.0}) is motivated by argument in \cite{Li}, and it is different from the proof in \cite{GL}.  
}
\end{remark}

\pf
{($\mathnormal{1}$). We first prove (\ref{2.11}). From the definition of $P_{N_0}(x, y, t)$, it is easy to see that $\lim_{t\rightarrow 0}P_{N_0}(x, y, t)= \delta_{y}(x)$. In particular, 
\begin{align}
F_{N_0}(x, y, t)&= H(x, y, t)- P_{N_0}(x, y, t) \nonumber \\
&= -\int_{0}^{t}\frac{\partial}{\partial s}\int_{M^n} H(x, z, t-s)P_{N_0}(z, y, s)d\mu(z)ds \nonumber\\
&= -\int_{0}^{t}\int_{M^n}\Big(\frac{\partial}{\partial s}- \Delta_{z}\Big)P_{N_0}(z, y, s)\cdot H(x, z, t- s)d\mu(z) ds\nonumber
\end{align}
where $\Delta_z$ is the Laplacian with respect to the $z$-variable.

From (\ref{2.7}) and the definition of $\eta$, when $z\in B(r)$,
\begin{equation}\label{2.10a}
{\Big|\Big(\frac{\partial}{\partial s}- \Delta_{z}\Big)P_{N_0}(z, y, s)\Big|\leq C_{1}s^3\exp{(-\frac{d^2(z, y)}{4s})}
}
\end{equation}
and when $z\in B(2r)\backslash B(r)$, 
\begin{equation}\label{2.10b}
{\Big|\Big(\frac{\partial}{\partial s}- \Delta_{z}\Big)P_{N_0}(z, y, s)\Big|\leq C_{2}s^{-\frac{n}{2}-1}\exp{(-\frac{d^2(z, y)}{4s})}
}
\end{equation}
Hence
\begin{align}
|F_{N_0}(x, y, t)|&\leq C_1\int_{0}^{t} s^3\int_{B(r)} H(x, z, t- s)\exp{\Big(-\frac{d^2(z, y)}{4s}\Big)} d\mu(z) ds\nonumber \\
& \quad +C_2\int_{0}^{t} s^{-\frac{n}{2}-1}\int_{B(2r)\backslash B(r)} H(x, z, t- s)\exp{\Big(-\frac{d^2(z, y)}{4s}\Big)} d\mu(z) ds \nonumber \\
&\leq (\mathit{a})+ (\mathit{b}) 
\end{align}
We can find $0< t_1\leq 1$ and $k_0> 0$, such that if $s\in (0, t_1)$, then 
\[V(B_{p}(\sqrt{s}))\geq k_0 s^{\frac{n}{2}} \quad for \ any \ p\in B_{y}(3r) \]
In the rest of the proof, assume $t\in (0, t_1]$, we have two cases.

Case (I): If $x\in B_{y}(3r)$ and $z\in B_{y}(2r)$, then from \cite{LY} and the above volume lower bound, 
\begin{align}
H(x, z, t- s)&\leq C V^{-\frac{1}{2}}\Big(B_x(\sqrt{t- s})\Big) V^{-\frac{1}{2}} \Big(B_z(\sqrt{t- s})\Big)\cdot \exp{\Big[CK(t- s)- \frac{6d^2(z, x)}{25(t- s)}\Big]}\label{2.12.0} \\
& \leq C(K, k_0, n)(t- s)^{-\frac{n}{2}} \exp{\Big(- \frac{6d^2(z, x)}{25(t- s)}\Big)}\label{2.12}
\end{align}

Case (II): If $x\notin B_y(3r)$ and $z\in B_y(2r)$, using (\ref{2.12.0}), $d(x, z)\geq r$ and volume comparison theorem, 
\begin{align}
H(x, z, t- s)&\leq C V^{-1}\Big(B_z(\sqrt{t- s})\Big) \cdot \Big[ \frac{V_{-K}\Big(\sqrt{t- s}+ d(x, z)\Big)}{V_{-K}\Big(\sqrt{t- s}\Big)} \Big]^{\frac{1}{2}} \nonumber \\
& \quad \cdot \exp{\Big[CK(t- s)- \frac{6d^2(z, x)}{25(t- s)}\Big]}\nonumber \\
& \leq C(K, k_0, n, r) \exp{\Big(- \frac{23d^2(z, x)}{100(t- s)}\Big)}\label{2.13}
\end{align}

Note in Case (I), $inj_{g}(x)$ has a uniform lower bound, hence it is easy to get
\begin{equation}\label{2.14}
{\int_{B_y(r)} s^{-\frac{n}{2}} \exp{\Big(-\frac{d^2(z, x)}{100s}\Big)} d\mu(z)\leq C
}
\end{equation}
for any $s\in (0, t_1]$.

Now using (\ref{2.12}), (\ref{2.13}), (\ref{2.14}) and the classical inequality 
\[\frac{d^2(x, z)}{t- s}+ \frac{d^2(y, z)}{s}\geq \frac{d^2(x, y)}{t}\]
we can get 
\begin{equation}\nonumber
{\int_{B_y(r)} H(x, z, t-s)\exp{\Big(-\frac{d^2(z, y)}{4s}\Big)} d\mu(z)\leq C\exp{\Big(-\frac{23d^2(x, y)}{100t}\Big)}
}
\end{equation}
Hence 
\begin{equation}\label{2.15}
{(\mathit{a})\leq C t^4\exp{\Big(-\frac{23d^2(x, y)}{100t}\Big)} 
}
\end{equation}

Similarly,
\begin{equation}\nonumber
{\int_{B_y(2r)\backslash B_y(r)} H(x, z, t-s)\exp{\Big(-\frac{d^2(z, y)}{4s}\Big)} d\mu(z)\leq C \exp{\Big(-\frac{3r^2}{100s}\Big)} \exp{\Big(-\frac{d^2(x, y)}{5t}\Big)}
}
\end{equation}
Hence 
\begin{align}
(\mathit{b})&\leq C_2\Big[\int_0^t s^{-\frac{n}{2}- 1} \exp{\Big(-\frac{3r^2}{100s}\Big)} ds\Big] \exp{\Big(-\frac{d^2(x, y)}{5t}\Big)} \nonumber\\
& \leq Ct^4 \exp{\Big(-\frac{d^2(x, y)}{5t}\Big)} \label{2.16}
\end{align}
By (\ref{2.15}) and (\ref{2.16}), (\ref{2.11}) is proved.

($\mathnormal{2}$). The strategy to prove (\ref{2.11.0}) is similar. 
\begin{align}
\frac{\partial}{\partial t}F_{N_0}(x, y, t)&= \frac{\partial}{\partial t} \Big[-\int_{0}^{t}\int_{M^n}\Big(\frac{\partial}{\partial s}- \Delta_{z}\Big)P_{N_0}(z, y, s)\cdot H(x, z, t- s)d\mu(z) ds \Big]\nonumber \\
&= -\int_{0}^{t}\int_{M^n}\Big(\frac{\partial}{\partial s}- \Delta_{z}\Big)P_{N_0}(z, y, s)\cdot \Big(\frac{\partial }{\partial t}H(x, z, t- s) \Big) d\mu(z) ds \nonumber \\
&\quad + \Big(\Delta_x- \frac{\partial}{\partial t}\Big)P_{N_0}(x, y, t) \nonumber \\
&= (\mathit{\gamma})+ (\mathit{\tau}) \nonumber
\end{align}

From (\ref{2.10a}), (\ref{2.10b}) and $P_{N_0}(x, y, t)= 0$ when $x\notin B(2r)$, 
\begin{equation}\label{2.17}
{(\mathit{\tau})\leq Ct^4\exp{\Big(-\frac{d^2(x, y)}{5t}\Big)}
}
\end{equation}

Now we estimate $(\mathit{\gamma})$.
\begin{align}
(\mathit{\gamma})&= -\int_{0}^{t}\int_{M^n}\Big(\frac{\partial}{\partial s}- \Delta_{z}\Big)P_{N_0}(z, y, s)\cdot \Big(\Delta_{z}H(x, z, t- s) \Big) d\mu(z) ds \nonumber \\
&= -\int_{0}^{t}\int_{M^n}\Big[\Delta_z \Big(\frac{\partial}{\partial s}- \Delta_{z}\Big)P_{N_0}(z, y, s)\Big] \cdot H(x, z, t- s)d\mu(z) ds \nonumber
\end{align}

Similar with (\ref{2.10a}) and (\ref{2.10b}), from (\ref{2.7}), when $z\in B(r)$, 
\begin{align}
\Big|\Delta_z \Big(\frac{\partial}{\partial s}- \Delta_{z}\Big)P_{N_0}(z, y, s)\Big|\leq C_3s\exp{\Big(-\frac{d^2(z, y)}{4s}\Big)} \label{2.10c}
\end{align}
and when $z\in B(2r)\backslash B(r)$, 
\begin{align}
\Big|\Delta_z \Big(\frac{\partial}{\partial s}- \Delta_{z}\Big)P_{N_0}(z, y, s)\Big|\leq C_4s^{-\frac{n}{2}- 3}\exp{\Big(-\frac{d^2(z, y)}{4s}\Big)} \label{2.10d}
\end{align}

Following similar argument in the proof of (\ref{2.11}), using (\ref{2.10c}), (\ref{2.10d}) instead of (\ref{2.10a}), (\ref{2.10b}), 
\begin{equation}\label{2.18}
{(\mathit{\gamma})\leq Ct^{2}\exp{\Big(-\frac{d^2(x, y)}{5t}\Big)} 
}
\end{equation}

From (\ref{2.17}) and (\ref{2.18}), 
\begin{equation}\nonumber
{\Big|\frac{\partial}{\partial t}F_{N_0}(x, y, t)\Big|\leq (\mathit{\gamma})+ (\mathit{\tau})\leq Ct^2\exp{(-\frac{d^2(x, y)}{5t})}
}
\end{equation}
}
\qed

\section{The short time asymptotics of $N(H, t)$}
From (\ref{2.6}) and (\ref{2.9}) in Theorem \ref{thm 2.1}, there exists $0< t_0\leq 1$ such that 
\begin{equation}\label{3.0}
{\frac{1}{2}\leq (4\pi t)^{\frac{n}{2}}\exp \Big(\frac{d^2(x, y)}{4t}\Big) H_{N_0}(x, y, t)\leq 2
}
\end{equation}
holds when $x\in B(\frac{r}{2})$ and $0< t\leq t_0$. In section $3$ and section $4$, we assume that $t\in (0, t_0]$ and $(M^n, g)$, $H$ are from Theorem \ref{thm 2.1}. 

\begin{lemma}\label{lem 3.1}
{\begin{equation}\label{3.1}
{\int_{B(\frac{r}{2})} \Big[\ln \frac{H(x, y, t)}{H_{N_0}(x, y, t)}\Big] \cdot H(x, y, t) d\mu(x)= O(t^{2})
}
\end{equation}
}
\end{lemma}

\pf
{Assume $x\in B(\frac{r}{2})$, $t\leq t_0$, then $P_{N_0}(x, y, t)= H_{N_0}(x, y, t)$. Hence 
\[F_{N_0}(x, y, t)= H(x, y, t)- H_{N_0}(x, y, t)\] 

From (\ref{2.11}), 
\begin{equation}\label{3.2}
{|F_{N_0}(x, y, t)|\leq Ct^{N_0+1- \frac{n}{2}} \exp \Big(-\frac{d^2(x, y)}{5t}\Big)
}
\end{equation}

If $F_{N_0}(x, y, t)> 0$, then 
\begin{align}
|\ln \frac{H}{H_{N_0}} \cdot H|(x, y, t) &= \ln \Big(1+ \frac{F_{N_0}}{H_{N_0}} \Big)\cdot H\leq \frac{F_{N_0}}{H_{N_0}}\cdot H \nonumber\\
& \leq Ct^{N_0+1} \exp \Big(\frac{d^2(x, y)}{20t}\Big)\cdot H(x, y, t) \nonumber
\end{align}
If $F_{N_0}(x, y, t)\leq 0$, then $H(x, y, t)\leq H_{N_0}(x, y, t)$,
\begin{align}
|\ln \frac{H}{H_{N_0}} \cdot H|(x, y, t) &= |\ln H(x, y, t)- \ln H_{N_0}(x, y, t)|\cdot H(x, y, t)\nonumber\\
&= |\frac{1}{\xi}[H(x, y, t)- H_{N_0}(x, y, t)]|\cdot H(x, y, t) \nonumber
\end{align}
where $H(x, y, t)\leq \xi\leq H_{N_0}(x, y, t)$. Hence 
\begin{equation}\nonumber
{|\ln \frac{H}{H_{N_0}} \cdot H|(x, y, t) \leq |\frac{F_{N_0}}{H}|\cdot H= F_{N_0}\leq Ct^{N_0+1- \frac{n}{2}} \exp \Big(-\frac{d^2(x, y)}{5t}\Big)
}
\end{equation}
By the above, 
\begin{equation}\label{3.3}
{|\ln \frac{H}{H_{N_0}} \cdot H|(x, y, t) \leq Ct^{4} \Big[t^{\frac{n}{2}}\exp \Big(\frac{d^2(x, y)}{20t}\Big)\cdot H+ \exp \Big(-\frac{d^2(x, y)}{5t}\Big) \Big]
}
\end{equation}
From (\ref{3.0}) and (\ref{3.2}), 
\begin{equation}\label{3.4}
{H(x, y, t)\leq |H_{N_0}|+ |F_{N_0}|\leq 2(4\pi t)^{-\frac{n}{2}}\exp \Big(-\frac{d^2(x, y)}{4t}\Big)+ Ct^{4}\cdot \exp \Big(-\frac{d^2(x, y)}{5t}\Big) 
}
\end{equation}
By (\ref{3.3}) and (\ref{3.4}), 
\begin{equation}\label{3.5}
{\Big|\ln \frac{H}{H_{N_0}} \cdot H\Big|(x, y, t) \leq Ct^{4}
}
\end{equation}
Hence $\int_{B(\frac{r}{2})} \Big[\ln \frac{H(x, y, t)}{H_{N_0}(x, y, t)}\Big] \cdot H(x, y, t) d\mu(x)= O(t^{2})$.
}
\qed

\bigskip

{\it \textbf{Proof of (\ref{3.6})}:}~
{\begin{align}
-\int_{M^n} fH d\mu = \int_{M^n\backslash B(\frac{r}{2})} (-f H)d\mu + \int_{B(\frac{r}{2})} (-f H) d\mu = (I)+ (II) \nonumber
\end{align}
Firstly, we estimate $(I)$. From \cite{LY}, we have
\begin{equation}\nonumber
{H(x, y, t)\leq C V^{-\frac{1}{2}}\Big(B_x(\sqrt{t})\Big) V^{-\frac{1}{2}} \Big(B_y(\sqrt{t})\Big)\cdot \exp{\Big[CKt- \frac{d^2(x, y)}{5t}\Big]}
}
\end{equation}
If $x\in M^n\backslash B(\frac{r}{2})$ and $t$ is small enough, using volume comparison theorem, 
\begin{equation}\label{3.5.1}
{H(x, y, t)\leq C V^{-1}\Big(B_y(\sqrt{t})\Big) \cdot \frac{V_{-K}\Big(\sqrt{t}+ d\Big)}{V_{-K}\Big(\sqrt{t}\Big)}\cdot \exp{\Big[CKt- \frac{d}{5t}\Big]}\leq Ct^2\exp{\Big(-\frac{d^2}{6t}\Big)}
}
\end{equation}
where $C$ depends on $n$, $K$, $r$ and the metric $g$ near $y$. Choose $t$ small enough such that $H\leq Ct^2\leq e^{-1}$, then by the monotonicity of $h(x)= \ln x \cdot x$ on $(0, e^{-1}]$, 
\begin{equation}\nonumber
{\Big|\ln H(x, y, t)\cdot H(x, y, t) \Big| \leq \Big|\ln \Big[Ct^2\exp{\Big(-\frac{d^2}{6t}\Big)}\Big]\cdot \Big[Ct^2\exp{\Big(-\frac{d^2}{6t}\Big)}\Big]\Big|
}
\end{equation}
hence
\begin{align}
|(I)|&= |\int_{M^n\backslash B(\frac{r}{2})} [\ln H+ \frac{n}{2}\ln (4\pi t)]\cdot H d\mu(x) | \nonumber \\
&\leq \int_{M^n\backslash B(\frac{r}{2})} \Big|\ln \Big[Ct^2\exp{\Big(-\frac{d^2}{6t}\Big)}\Big]\cdot \Big[Ct^2\exp{\Big(-\frac{d^2}{6t}\Big)}\Big]\Big| d\mu(x) \nonumber \\
&\quad + \frac{n}{2} \int_{M^n\backslash B(\frac{r}{2})} |\ln (4\pi t)\cdot \Big[Ct^2\exp{\Big(-\frac{d^2}{6t}\Big)}\Big]| d\mu(x) \nonumber \\
&\leq O(t^{\frac{3}{2}}) \label{use}
\end{align}
in the last inequality, we used $Rc\geq -Kg$ and volume comparison theorem.

\begin{align}
|(II)|&= \int_{B(\frac{r}{2})} [\ln H+ \frac{n}{2}\ln (4\pi t)]\cdot H d\mu\nonumber \\
&= \int_{B(\frac{r}{2})} \ln \frac{H}{H_{N_0}}\cdot H d\mu(x)+ \int_{B(\frac{r}{2})} \Big[\ln H_{N_0}+ \frac{n}{2}\ln (4\pi t)\Big]\cdot H d\mu(x) \nonumber \\
&= (III)+ (IV) \nonumber
\end{align}
By Lemma \ref{lem 3.1}, $(III)= O(t^2)$. From Lemma \ref{lem 3.3} in the following, 
\[(IV)= -\frac{n}{2}+ \frac{1}{2}R(y)\cdot t + O(t^{\frac{3}{2}})\]
By all the above, we get our conclusion.
}
\qed

\begin{lemma}\label{lem 3.3}
{\begin{equation}\label{3.9}
{\int_{B(\frac{r}{2})} \Big[\ln H_{N_0}+ \frac{n}{2}\ln (4\pi t)\Big]\cdot H d\mu(x)= -\frac{n}{2}+ \frac{1}{2}R(y)\cdot t + O(t^{\frac{3}{2}})
}
\end{equation}
}
\end{lemma}

\pf
{$(I)\doteqdot \int_{B(\frac{r}{2})} \Big[\ln H_{N_0}+ \frac{n}{2}\ln (4\pi t)\Big]\cdot H d\mu(x)$ in the following proof. From Theorem \ref{thm 2.1}, 
\begin{equation}\nonumber
{\ln H_{N_0}= -\frac{n}{2}\ln (4\pi t)- \frac{d^2(x, y)}{4t}+ \ln \Big(\sum_{k= 0}^{N_0}\varphi_k t^k\Big)
}
\end{equation}
and 
\begin{equation}\nonumber
{\ln \Big(\sum_{k= 0}^{N_0}\varphi_k t^k\Big)= \ln \varphi_0+ \frac{\varphi_1}{\varphi_0}\cdot t+ O(t^2)
}
\end{equation}
Hence 
\begin{equation}\nonumber
{(I)= \int_{B(\frac{r}{2})}\Big[-\frac{d^2(x, y)}{4t}+ \ln \varphi_0+ \frac{\varphi_1}{\varphi_0}\cdot t+ O(t^2)\Big]\cdot H d\mu(x)
}
\end{equation}
Now using $(iii)$ of Theorem \ref{thm 2.1}, 
\begin{align}
(I)&= \int_{B(\frac{r}{2})} \Big[-\frac{d^2}{4t}+ \frac{1}{12}R_{pq}(y)x^px^q+ O(d^3)+ \Big(\frac{R(y)}{6}+ O(d)\Big)t\Big]\cdot H d\mu(x)+ O(t^2) \nonumber \\
&= (II)+ (III)+ (IV)+ (V)+ (VI)+ O(t^2) \nonumber 
\end{align}
where 
\begin{align}
(II)&= \int_{B(\frac{r}{2})} \Big(-\frac{d^2(x, y)}{4t}\Big)\cdot H d\mu(x), \quad (III)= \frac{1}{12}\int_{B(\frac{r}{2})} \Big(R_{pq}(y)x^px^q)\cdot H d\mu(x), \nonumber \\
(IV)&= C\int_{B(\frac{r}{2})} d^3(x, y)\cdot H(x, y, t) d\mu(x), \quad (V)= \frac{R(y)}{6}t\cdot \int_{B(\frac{r}{2})} H(x, y, t) d\mu(x), \nonumber \\
(VI)&= Ct\cdot \int_{B(\frac{r}{2})} d(x, y)\cdot H(x, y, t) d\mu(x) \nonumber
\end{align}

From (\ref{3.5.1}), 
\begin{equation}\nonumber
{\int_{B(\frac{r}{2})} H= \int_{M^n} H- \int_{M^n\backslash B(\frac{r}{2})} H= 1+ O(t^2)
}
\end{equation}
Hence 
\begin{equation}\nonumber
{(V)= \frac{1}{6}R(y)\cdot t+ O(t^2)
}
\end{equation}
Using (\ref{3.4}) and the fact that
\begin{equation}\nonumber
{\int_{R^n}O(|x|^{k})(4\pi t)^{-\frac{n}{2}} \exp (-\frac{|x|^2}{4t}) dx= O(t^{\frac{k}{2}})
}
\end{equation}
where $k$ is any nonnegative integer, we can get $(IV)= O(t^{\frac{3}{2}})$ and $(VI)= O(t^{\frac{3}{2}})$. Similarly, 
\begin{equation}\nonumber
{(III)= \frac{1}{6}R(y)\cdot t+ O(t^2)
}
\end{equation}
Finally, from the following Lemma \ref{lem 3.4}, 
\begin{equation}\nonumber
{(II)= -\frac{n}{2}+ \frac{1}{6}R(y)\cdot t+ O(t^{\frac{3}{2}})
}
\end{equation}
By all the above, the conclusion is proved.
}
\qed

\begin{lemma}\label{lem 3.4}
{\begin{equation}\label{3.10}
{-\frac{1}{4t}\int_{B(\frac{r}{2})} d^2(x, y)\cdot H d\mu(x)= -\frac{n}{2}+ \frac{1}{6}R(y)\cdot t+ O(t^{\frac{3}{2}})
}
\end{equation}
}
\end{lemma}

\pf
{$(II)\doteqdot -\frac{1}{4t}\int_{B(\frac{r}{2})} d^2(x, y)\cdot H d\mu(x)$, then
\begin{equation}\label{3.11}
{(II)= -\frac{1}{4t}\int_{B(\frac{r}{2})} d^2(x, y)\cdot (H_{N_0}+ F_{N_0})\cdot \alpha dx
}
\end{equation}
where $dx$ in the integral of (\ref{3.11}) is the volume element of Euclidean space $\mathbb{R}^n$, and
\begin{equation}\nonumber
{\alpha= \sqrt{det(g)}= 1- \frac{1}{6}R_{pq}(y)x^px^q+ O(d^3(x, y))
}
\end{equation}
Then
\begin{align}
(II)&= -\frac{1}{4t}\int_{B(\frac{r}{2})} d^2(x, y)\cdot (4\pi t)^{-\frac{n}{2}} \exp{\Big(-\frac{d^2(x, y)}{4t}\Big)} (\varphi_0+ \varphi_1 t)\cdot \alpha dx+ O(t^2) \nonumber \\
&= \Big[-\frac{1}{4t}- \frac{1}{24}R(y) \Big]\cdot \int_{B(\frac{r}{2})}d^2\cdot (4\pi t)^{-\frac{n}{2}}\exp (-\frac{d^2}{4t}) dx \nonumber \\
&\quad + \frac{1}{48t}\int_{B(\frac{r}{2})} (4\pi t)^{-\frac{n}{2}} (R_{pq}(y) x^px^q) d^2\cdot \exp (-\frac{d^2}{4t}) dx+ O(t^{\frac{3}{2}}) \nonumber \\
&= \Big[-\frac{1}{4t}- \frac{1}{24}R(y)\Big]\cdot 2nt+ \frac{1}{48t}I_n+ O(t^{\frac{3}{2}}) \nonumber
\end{align}
where 
\begin{equation}\nonumber
{I_n= \int_{\mathbb{R}^n} (4\pi t)^{-\frac{n}{2}}\Big(\sum_{k=1}^{n}\lambda_k x_k^2\Big)\cdot \Big(\sum_{i= 1}^{n}x_i^2\Big) \exp \Big(-\frac{1}{4t}\cdot \sum_{j= 1}^n x_j^2\Big) dx
}
\end{equation}
in above we diagonalize $R_{pq}(y)$ and let $\lambda_k= R_{kk}(y)$.

We can get $I_1= 12\lambda_1 t^2$, and the induction formula 
\[I_n= I_{n- 1}+ 4(\sum_{i=1}^{n}\lambda_i)t^2+ 4(n+ 1)\lambda_n t^2 \]
Then it is easy to get 
\begin{equation}\label{I_n}
{I_n= 4(n+ 2)(\sum_{i=1}^{n}\lambda_i)t^2= 4(n+ 2)R(y)t^2
}
\end{equation}
By all the above $(II)= -\frac{n}{2}+ \frac{R(y)}{6}t+ O(t^{\frac{3}{2}})$, the lemma is proved.
}
\qed

\section{The short time asymptotics of $\frac{\partial}{\partial t}\Big[N(H, t) \Big]$}

To study $\frac{\partial}{\partial t}\Big[N(H, t) \Big]$, we need to switch the differentiation with integration firstly. Because the manifold $M^n$ can be non-compact, we need to be more careful on the switch  of the order of differentiation and integration. The next lemma justifies this switch in our case.

\begin{lemma}\label{switch}
{\begin{align}
\frac{\partial}{\partial t} \Big[\int_{M^n} H(-f) d\mu(x)\Big]= \int_{M^n} \Big[H(-f)\Big]_{t} d\mu(x)
\end{align}
}
\end{lemma}

\pf
{Define $\varphi_{\rho}(x)= \phi\Big(\frac{d(x, y)}{\rho}\Big)$, where $\phi$ is defined in Appendix, $\rho> 1$ is a constant. Fix $t> 0$, define $G(x, t)= [H(-f)](x, y, t)$. For any $\epsilon> 0$, assume $1> l> 0$ (if $l< 0$, similar argument works). Then
\begin{align}
&\Big|\int_{M^n} \frac{G(x, t+ l)- G(x, t)}{l} d\mu(x)- \int_{M^n} G_{t}\varphi_{\rho} d\mu(x)\Big| \nonumber \\
&\leq \int_{B(\rho)} |G_{t}(x, t+ \xi_{x}l)- G_{t}(x, t)| d\mu(x)+ 2\int_{M^n\backslash B(\rho)} \sup_{s\in [t, t+ l]} |G_{t}(x, s)| d\mu(x) \nonumber \\
&\leq \int_{B(\rho)} \Big|\frac{\partial^2}{\partial^2 t}G(x, t+ \zeta_{x}l) \Big| d\mu(x)\cdot l+ 2\int_{M^n\backslash B(\rho)} \Big(\sup_{s\in [t, t+ l]} |G_{t}(x, s)| \Big)d\mu(x) \nonumber \\
&\leq (I)+ (II) 
\end{align}

We firstly estimate $(II)$. From \cite{LY}, for $s\in [t, t+ l]$, 
\begin{align}
\frac{H_{t}}{H}(x, y, s)\geq \frac{1}{2}\Big[\frac{|\nabla H|^2}{H^2}- \frac{2n}{s}- CK\Big]\geq -\frac{C}{s} \label{4.0.2}
\end{align}
where $C= C(K, n)$. From Corollary \ref{cor 2.5}, 
\begin{align}
\frac{H_{t}}{H}(x, s)&\leq \frac{2}{s}\Big\{n+ (4+ Ks)\ln \Big[\frac{C(K, t+ 1)}{H(x, y, s) V^{\frac{1}{2}}(B_{x}(\sqrt{\frac{s}{2}})) V^{\frac{1}{2}}(B_{y}(\sqrt{\frac{s}{2}}))}\Big]\Big\} \nonumber \\
&\leq \frac{C}{s}\Big(1+ |\ln H|+ \Big|\ln \Big[V\Big(B_{x}\Big(\sqrt{\frac{s}{2}}\Big)\Big)\cdot V\Big(B_{y}\Big(\sqrt{\frac{s}{2}}\Big)\Big)\Big]\Big|\Big) \label{4.0.3}
\end{align}

When $x\in M^n\backslash B(\rho)$, using volume comparison theorem, 
\begin{align}
\Big|\ln \Big[V\Big(B_{x}\Big(\sqrt{\frac{s}{2}}\Big)\Big)\cdot V\Big(B_{y}\Big(\sqrt{\frac{s}{2}}\Big)\Big)\Big]\Big| &\leq 2|\ln V\Big(B_{y}\Big(\sqrt{\frac{s}{2}}\Big)\Big) | \nonumber \\
& \quad + \Big|\ln V_{-K}\Big(\sqrt{\frac{s}{2}}\Big) \Big| + |\ln V_{-K}(s+ d(x, y))| \nonumber \\
&\leq C(|\ln s|+ s+ d) \label{4.0.4}
\end{align}

where $C$ is independent of $\rho$. From (\ref{4.0.2}), (\ref{4.0.3}) and (\ref{4.0.4}), 
\begin{align}
\Big|\frac{H_t}{H}\Big|(x, s)\leq \frac{C}{s} (|\ln H|+ |\ln s|+ s+ d) \label{4.0.5}
\end{align}
when $x\in M^n\backslash B(\rho)$.

From (\ref{4.0.5}), on $M^n\backslash B(\rho)$,
\begin{align}
|(-f)H_{t}|(x, s) &\leq  \Big[|\ln H|+ \frac{n}{2}|\ln (4\pi s)|\Big] \cdot |H_{t}|(x, s) \nonumber \\
& \leq \Big[|\ln H|+ \frac{n}{2}|\ln (4\pi s)|\Big] \cdot C|H|\cdot s^{-1}(|\ln H|+ |\ln s|+ s+ d) \nonumber \\
&\leq \frac{C}{s} \cdot H\Big[|\ln H|^2+ |\ln s|^2+ s^2+ d^2\Big]  \label{4.0.5.1}
\end{align}

From (\ref{4.0.5}) and (\ref{4.0.5.1}), if $s\in [t, t+ l]$ and $x\in M^n\backslash B(\rho)$, 
\begin{align}
|G_{t}(x, s)|&\leq \Big[|H_{t}|+ \frac{n}{2s}|H|+ |(-f)H_{t}| \Big](x, s) \nonumber \\
&\leq \frac{C}{s}H\cdot \Big(|\ln H|+ |\ln s|+ s+ d\Big)\nonumber \\
&\quad + \frac{n}{2s}|H|+ \frac{C}{s}H\cdot \Big(|\ln H|^2+ |\ln s|^2+ s^2+ d^2\Big) \nonumber \\
&\leq \frac{C}{s}H\cdot \Big(|\ln H|^2+ |\ln s|^2+ s^2+ d^2\Big) \label{s2}
\end{align}
where $C$ is independent of $\rho$. We can choose $l$ smooth enough such that $(t+ l)\leq 2t$, then using (\ref{3.5.1}) and (\ref{s2}), on $x\in M^n\backslash B(\rho)$ 
\begin{align}
|G_{t}(x, s)|&\leq Cs\exp{\Big(-\frac{d^2}{6s}\Big)}\cdot \Big[\Big|C+ 2\ln s- \frac{d^2}{6s}\Big|^2+ |\ln s|^2+ s^2+ d^2\Big] \nonumber \\
&\leq C(t+ l)\exp{\Big(-\frac{d^2}{6(t+ l)}\Big)}\cdot \Big[t^2+ d^2+ |\ln t|^2+ \Big(\frac{d^2}{t}\Big)^2\Big] \nonumber \\
&\leq Ct\exp{\Big(-\frac{d^2}{12t}\Big)}\cdot \Big[t^2+ |\ln t|^2+ \Big(\frac{d^2}{t}\Big)^2\Big] 
\end{align}

Hence for any $\epsilon> 0$, we can find $\rho_0> 1$, such that if $\rho\geq \rho_0$,
\begin{align}
\int_{M^n\backslash B(\rho)} \Big(\sup_{s\in [t, t+ l]} |G_{t}(x, s)| \Big)d\mu(x)< \frac{\epsilon}{4} \label{s3}
\end{align}

On the other hand, because $0< l< 1$,
\begin{align}
\int_{B(\rho)} |G_{tt}(x, t+ \zeta_x l)|d\mu(x) &\leq \int_{B(\rho)} \sup_{s\in [t, t+ 1]}|G_{tt}(x, s)|d\mu(x) \nonumber \\
&\leq C(\rho) \label{s4}
\end{align}

Choose $l\leq \frac{\epsilon}{4C(\rho)}$, from (\ref{s3}) and (\ref{s4}), if $\rho> \rho_0$,
\begin{align}
\Big|\int_{M^n} \frac{G(x, t+ l)- G(x, t)}{l} d\mu(x)- \int_{M^n} G_{t}\varphi_{\rho} d\mu(x)\Big|< \epsilon \label{s5}
\end{align}

It is easy to see from Lemma \ref{lem 4.0} and its proof, $\lim_{\rho\rightarrow \infty}\int_{M^n} G_{t}\phi_{\rho}$ exists and 
\begin{align}
\lim_{\rho\rightarrow \infty}\int_{M^n} G_{t}\phi_{\rho}= \int_{M^n} G_{t}\label{s6}
\end{align}

From (\ref{s5}) and (\ref{s6}), we get our conclusion.
}
\qed

By Cheng-Li-Yau's result (see \cite{CLY}), $\lim_{t\rightarrow 0}t\ln H= -\frac{d^2}{4}$ and the limit is uniform for any $x$ in $B(r)$. Hence we can assume
\begin{equation}\nonumber
{t\ln H(x, y, t)= -\frac{d^2(x, y)}{4}+ \epsilon (t, x, y)
}
\end{equation}
Sometimes $\epsilon(t, x, y)$ will be simplified as $\epsilon$, then
\begin{equation}\label{4.02}
{t(-f)= \frac{n}{2}t\ln (4\pi t)- \frac{d^2}{4}+ \epsilon
}
\end{equation}
where $\lim_{t\rightarrow 0}\epsilon(t, x, y)= 0$, and the limit is uniform for any $x$ in $B(r)$. Without losing generality, we can assume that $\varphi_{0}(x, y)\geq \frac{1}{2}$ when $x\in B(\frac{r}{2})$.

\begin{lemma}\label{lem 4.1}
{\begin{equation}\label{4.1}
{\int_{B(\frac{r}{2})} E(-f) d\mu(x)= -\frac{n}{2}+ \frac{1}{3}R(y)\cdot t+ o(t)
}
\end{equation}
and 
\begin{equation}\label{4.2}
{\int_{B(\frac{r}{2})} E(-f) O(d(x, y)) d\mu(x)= o(1)
}
\end{equation}
where $\lim_{t\rightarrow 0}\frac{o(t)}{t}= 0$.
}
\end{lemma}

\pf
{\begin{align}
&\int_{B} E(-f)d\mu(x)= \int_{B} \frac{H_{N_0}}{\sum_{k= 0}^{N_{0}}\varphi_k t^k}\cdot (-f) d \mu(x) \nonumber \\
&\quad = \int_B \Big(\frac{1}{\varphi_0}- \frac{\varphi_1}{\varphi_0^2}t\Big)H(-f) d\mu(x) + o(t) \nonumber \\
&\quad = \int_{B} \Big(1+ \frac{1}{12}R_{pq}(y)x^px^q- \frac{R(y)}{6}t\Big)H(-f)+ o(t) \nonumber \\
&\quad = -\frac{n}{2}+ \Big(\frac{1}{2}+ \frac{n}{12}\Big)R(y)t+ \frac{1}{12}\int_{B} R_{pq}(y)x^px^q\cdot H(-f) d\mu(x)+ o(t) \label{4.1.0}
\end{align}
in the last equation, we used (\ref{3.6}).

We estimate the third term on the right side of (\ref{4.1.0}).
\begin{align}
(I)&\doteqdot \frac{1}{12}\int_{B} R_{pq}(y)x^px^q\cdot H(-f) d\mu(x) \nonumber \\
& = \frac{1}{12} \int_{B}  R_{pq}(y)x^px^q\cdot H\Big[\ln H_{N_0}+ \frac{n}{2}\ln (4\pi t) \Big] d\mu(x) \nonumber \\
& =  \frac{1}{12} \int_{B} R_{pq}(y) x^px^q\cdot H_{N_0} \Big[-\frac{d^2}{4t}+ \ln \varphi_0 \Big]\cdot \alpha dx+ o(t) \nonumber \\
& = -\frac{1}{48t} \int_{B} E\cdot d^2\cdot R_{pq}(y) x^px^q dx+ o(t) \nonumber \\ 
& = -\frac{n+ 2}{12}R(y)t+ o(t)     \label{4.1.1}                                                  
\end{align}              

In the last equation above, we used (\ref{I_n}). From (\ref{4.1.0}) and (\ref{4.1.1}), we get (\ref{4.1}). To prove (\ref{4.2}), we will follow similar strategy.
\begin{align}
&\int_{B} E(-f) O(d) d\mu(x)= \int_B \Big(\frac{1}{\varphi_0}- \frac{\varphi_1}{\varphi_0^2}t\Big)H(-f) O(d) d\mu(x) + o(1) \nonumber \\
&\quad = \int_{B} H_{N_0}\Big[\ln H_{N_0}+ \frac{n}{2}\ln (4\pi t)\Big] O(d) d\mu(x)+ o(1) \nonumber \\
&\quad = \int_{B} E\Big(-\frac{d^2}{4t}+ \ln \varphi_0 \Big) O(d) d\mu(x)+ o(1)= o(1) \nonumber
\end{align}
(\ref{4.2}) is proved.
}
\qed

\begin{lemma}\label{lem 4.0}
{\begin{align}
\int_{M^n\backslash B} |(-f)H_{t}| d\mu(x)= O(t^{\frac{1}{2}}) \nonumber
\end{align}
where $t< < 1$ is small enough.
}
\end{lemma}

\pf
{Similarly as (\ref{4.0.5.1}), on $M^n\backslash B$, 
\begin{align}
|(-f)H_{t}|\leq \frac{C}{t} \cdot H\Big[|\ln H|^2+ |\ln t|^2+ t^2+ d^2\Big]  \nonumber
\end{align}

Hence 
\begin{align}
\int_{M^n\backslash B} |(-f)H_{t}|&\leq \frac{C}{t} \int_{M^n\backslash B} H\cdot |\ln H|^2+ \frac{C}{t} \int_{M^n\backslash B} H (|\ln t|^2+ t^2+ d^2)  \nonumber \\
& =(I)+ (II) \nonumber
\end{align}

Using (\ref{3.5.1}), volume comparison theorem and monotonicity of $h(x)= x(\ln x)^2$ when $x\in (0, e^{-2}]$, similar to the proof of (\ref{use}), 
\begin{align}
(I)\leq O(t^{\frac{1}{2}}) \nonumber
\end{align}

Using (\ref{2.11}), when $x\in M^n\backslash B$, 
\begin{equation}\label{4.2.1}
{H\leq |\eta H_{N_0}|+ |F_{N_0}|\leq C\Big[t^{-\frac{n}{2}}\exp \Big(-\frac{d^2}{4t}\Big)+ t^{4}\cdot \exp \Big(-\frac{d^2}{5t}\Big) \Big]= O(t^2)\tilde{E}
}
\end{equation}

From (\ref{4.2.1}), it is easy to get
\begin{align}
(II)\leq O(t) \nonumber
\end{align}

By all the above, we get our conclusion.
}
\qed

\bigskip

{\it \textbf{Proof of (\ref{3.6.1})}:}~
{\begin{align}
\frac{\partial}{\partial t}\Big[\int_{M^n} H(-f) d\mu(x)\Big]&= \int_{M^n} \Big[H_t+ \frac{n}{2t}H+ (-f)H_t\Big]d\mu(x) \nonumber \\
&= \frac{n}{2t}+ \int_{M^n\backslash B(\frac{r}{2})} (-f)H_t d\mu(x)+ \int_{B(\frac{r}{2})} (-f)H_t d\mu(x) \nonumber \\
&= \frac{n}{2t}+ (I)+ (II) \nonumber
\end{align}

From Lemma \ref{lem 4.0} in the above, we have
\begin{align}
(I)= O(t^{\frac{1}{2}}) \nonumber
\end{align}

From Lemma \ref{lem 4.3} in the following, we get 
\begin{align}
(II)= -\frac{n}{2t}+ \frac{1}{2}R(y)+ o(1) \nonumber
\end{align}

From all the above, we get (\ref{3.6.1}).
}
\qed

\begin{lemma}\label{lem 4.3}
{\begin{equation}\nonumber
{\int_{B} (-f)H_t d\mu(x)= -\frac{n}{2t}+ \frac{1}{2}R(y)+ o(1)
}
\end{equation}
}
\end{lemma}

\pf
{From (\ref{2.11.0}) and (\ref{4.02}),
\begin{equation}\nonumber
{\int_{B} (-f)H_t d\mu(x)= \int_{B} (-f)\cdot (H_{N_0})_t+ O(t)
}
\end{equation}

\begin{align}
&\int_{B} (-f)(H_{N_0})_t d\mu(x)= \int_{B} \Big(\frac{d^2}{4t^2}- \frac{n}{2t}\Big) H_{N_0}\cdot (-f)d\mu(x)+ \int_{B} E\varphi_1 (-f) d\mu(x)+ o(1) \nonumber \\
&\quad = \frac{1}{4t^2}\int_{B} H_{N_0}(-f)d^2 d\mu(x)- \frac{n}{2t}\int_{B}H_{N_0}(-f)d\mu(x)+ \int_{B} E\varphi_1(-f) d\mu(x) + o(1)\nonumber \\
&\quad = (I)+ (II)+ (III)+ o(1)  \nonumber
\end{align}

Using Lemma \ref{lem 4.1}, 
\begin{align}
(III)&= \int_{B}E\varphi_1 (-f) d\mu(x)= \frac{1}{6}R(y)\int_{B} E(-f) d\mu(x)+ \int_{B} E(-f)\cdot O(d)\nonumber \\
&= -\frac{n}{12}R(y)+ o(1)  \nonumber
\end{align}

From (\ref{2.11}) and (\ref{3.6}),
\begin{align}
(II)&= -\frac{n}{2t}\int_{B} H_{N_0}(-f) d\mu(x) = -\frac{n}{2t}\int_{B} H(-f)d\mu(x) -\frac{n}{2t}\int_{B} O(t^{N_0+ 1})\tilde{E}(-f) d\mu(x) \nonumber \\
&= \frac{n^2}{4t}- \frac{n}{4}R(y)+ o(1)  \nonumber
\end{align}

Similarly, by the following Lemma \ref{lem 4.4}, 
\begin{align}
(I)&= \frac{1}{4t^2} \int_{B} (H+ O(t^{N_0+ 1})\tilde{E})(-f)\cdot d^2 d\mu(x)= \frac{1}{4t^2} \int_{B} H(-f)\cdot d^2 d\mu(x)+ o(1)\nonumber \\
&= -\frac{n(n+ 2)}{4t}+ \Big(\frac{n}{3}+ \frac{1}{2}\Big)R(y)+ o(1)  \nonumber
\end{align}

From all the above,
\begin{equation}\nonumber
{\int_{B} (-f)H_t d\mu(x)= -\frac{n}{2t}+ \frac{1}{2}R(y)+ o(1)
}
\end{equation}
}
\qed

\begin{lemma}\label{lem 4.4}
{\begin{equation}\nonumber
{\frac{1}{4t^2} \int_{B} H(-f)\cdot d^2 d\mu(x)= -\frac{n(n+ 2)}{4t}+ \Big(\frac{n}{3}+ \frac{1}{2}\Big) R(y)+ o(1)
}
\end{equation}
}
\end{lemma}

\pf
{We will use the similar strategy as in the proof of (\ref{3.6}). 
\begin{align}
\frac{1}{4t^2} \int_{B} H(-f)\cdot d^2 d\mu(x)&= \frac{1}{4t^2} \int_{B} \Big[\ln H_{N_0}+ \frac{n}{2}\ln (4\pi t)\Big]Hd^2 d\mu(x) \nonumber \\
&\quad + \frac{1}{4t^2} \int_{B} \Big[\ln \frac{H}{H_{N_0}}\Big]Hd^2 d\mu(x) \nonumber
\end{align}

From (\ref{3.5}), 
\begin{equation}\nonumber
{\Big[\ln \frac{H}{H_{N_0}}\Big]H= O(t^4)
}
\end{equation}

Hence, 
\begin{align}
&\frac{1}{4t^2} \int_{B} H(-f)\cdot d^2 d\mu(x)= \frac{1}{4t^2} \int_{B} \Big[-\frac{d^2}{4t}+ \ln \varphi_0+ \frac{\varphi_1}{\varphi_0}t+ O(t^2)\Big]Hd^2\cdot \alpha dx+ o(1) \nonumber\\
&= \frac{1}{4t^2}\int_{B} \Big(-\frac{d^2}{4t}+ \frac{1}{12}R_{pq}(y)x^px^q+ \frac{1}{6}R(y)t- \frac{R(y)}{24}d^2+ \frac{1}{48t}R_{pq}(y)x^px^q\cdot d^2\Big) \nonumber \\
&\quad \quad \quad\cdot E d^2 dx+ o(1) \nonumber \\
&\quad = -\frac{n(n+ 2)}{4t}+ \frac{-n^2+ 2n+ 4}{24}R(y)+ \frac{1}{192t^3}\int_{\mathbb{R}^n} ER_{pq}(y)x^px^q \cdot d^4 dx+ o(1) \nonumber
\end{align}

Define
\begin{equation}\nonumber
{Q_{n}= \int_{\mathbb{R}^n} ER_{pq}(y)x^px^q \cdot d^4 dx= \int_{\mathbb{R}^n} E\cdot \Big(\sum_{i= 1}^n \lambda_i x_i^2\Big)\cdot (\sum_{j= 1}^n x_j^2) dx  
}
\end{equation}
where we diagonalize $R_{pq}(y)$ and let $\lambda_i= R_{ii}(y)$. We can get $Q_1= 120\lambda_1 t^3$ and the induction formula:
\begin{equation} \nonumber
{Q_n= Q_{n- 1}+ 8(2n+ 5)\Big(\sum_{i= 1}^n \lambda_i \Big)t^3+ 8(n^2+ 4n+ 3)\lambda_n \cdot t^3
}
\end{equation}
Then it is easy to get $Q_n= 8(n^2+ 6n+ 8)R(y)\cdot t^3$, hence
\begin{equation}\nonumber
{\frac{1}{4t^2} \int_{B} H(-f)\cdot d^2 d\mu(x)= -\frac{n(n+ 2)}{4t}+ \Big(\frac{n}{3}+ \frac{1}{2}\Big) R(y)+ o(1)
}
\end{equation}
}
\qed

\appendix
\section{}
In \cite{Ham}, Richard Hamilton established an upper bound of Laplacian of positive solution to the heat equation on closed manifolds. We will generalize his theorem to complete manifolds with Ricci curvature bounded below. Our proof follows the similar strategy as \cite{Kots}. We firstly establish a preliminary estimate on $t|\Delta u|$ so that the maximum principle of Ni and Tam \cite{NT} may be applied to the quantity of interest in Hamilton's second derivative estimate. 

We introduce a cut-off function $\phi$ defined on $\mathbb{R}$, which is a smooth nonnegative nonincreasing funciton, is $1$ on $(-\infty, 1)$ and $0$ on $[2, +\infty)$. We can further assume
\begin{equation}\label{2.19}
{|\phi '|\leq 2 \ , \quad \quad  |\phi ''|+ \frac{(\phi ')^2}{\phi}\leq 16
}
\end{equation}

To prove the following Bernstein-type local estimate, we employ a technique of W.-X. Shi \cite{Shi} from the estimation of derivatives of curvature under the Ricci flow (see also \cite{RFAA}), define $F= (C+ t|\nabla u|^2)t^2|\Delta u|^2$, and consider the evolution of $F$.

\begin{lemma}\label{lem 2.3}
{Suppose $(M^n, g)$ is a complete Riemannian manifold. If $|u(x, t)|\leq \mathscr{M}$ is a solution to the heat equation on $B_{p}(4\rho)\times [0, T]$ for some $p\in M^n$, constants $\mathscr{M}$, $\rho$, $T$, $K> 0$, and $Rc\geq -Kg$ on $B_{p}(4\rho)$. Then 
\begin{equation}\label{2.20}
{t|\Delta u|\leq C(n, K, \mathscr{M})\Big[1+ T\Big(1+ \frac{1}{\rho^2}\Big)\Big] \cdot \Big(\frac{1}{\rho}+ 1\Big) \cdot \Big[T+ \coth \Big(\sqrt{\frac{K}{n- 1}}\rho \Big)\Big]
}
\end{equation}
holds on $B_{p}(\rho)\times [0, T]$.
}
\end{lemma}

\pf
{From \cite{Kots}, we get 
\begin{equation}\label{2.21}
{t|\nabla u|^2\leq C_1\Big[1+ T(1+ \frac{1}{\rho^2})\Big]\doteqdot C_2
}
\end{equation}
holds on $B_{p}(2\rho)\times [0, T]$, where $C_1= C_1(K, \mathscr{M})$. Define $C_3= 4C_2$, and 
\begin{equation}\nonumber
{F(x, t)= (C_3+ t|\nabla u(x, t)|^2)t^2 |\Delta u(x, t)|^2
}
\end{equation}
Long but straightforward computation gives
\begin{align}
\Big(\frac{\partial}{\partial t}- \Delta \Big)F&= -2(C_3+ t|\nabla u|^2)|\nabla \Delta u|^2- 8t^3 \sum_{i, j}\nabla_{i}\nabla_j u\nabla_i \Delta u\nabla_j u \Delta u \nonumber \\
&\quad -2t^3|\nabla^2 u|^2\cdot |\Delta u|^2+ 2t(C_3+ t|\nabla u|^2) |\Delta u|^2 \nonumber \\
&\quad + \Big[|\nabla u|^2- 2t Rc(\nabla u, \nabla u)\Big]t^2|\Delta u|^2 \nonumber
\end{align}
When $x\in B_{p}(4\rho)$, using $t|\nabla u|^2\leq C_2= \frac{1}{4}C_3$ and $Rc\geq -Kg$, 
\begin{align}
\Big(\frac{\partial}{\partial t}- \Delta \Big)F&\leq -10t^3|\nabla u|^2\cdot|\nabla \Delta u|^2+ 8t^3 |\nabla u| \cdot |\nabla \Delta u| \cdot |\nabla^2 u|\cdot |\Delta u| \nonumber \\
&\quad -2t^3|\nabla^2 u|^2\cdot |\Delta u|^2+ C_4t|\Delta u|^2 \nonumber \\
&\leq -\frac{2}{5}t^3|\nabla^2 u|^2\cdot |\Delta u|^2+ C_4t|\Delta u|^2 \nonumber
\end{align}
where $C_4= (2KT+ 11)C_2$, the term with coefficient $-\frac{2}{5}$ arose from the inequality $-10x^2+ 8xy- 2y^2\leq -\frac{2}{5}y^2$. On the other hand, we know $|\nabla^2 u|^2\geq \frac{1}{n}|\Delta u|^2$, hence
\begin{align}
\Big(\frac{\partial}{\partial t}- \Delta \Big)F&\leq -\frac{2}{5n}t^3|\Delta u|^4+ C_4t|\Delta u|^2 \nonumber \\
&\leq -\frac{1}{5nt}\Big[t^2|\Delta u|^2\Big]^2+ \frac{5n}{4t}C_4 \nonumber \\
&\leq -\frac{C_6}{t}F^2+ \frac{C_5}{t} \nonumber
\end{align}
in the last equality we used $F\leq (C_3+ C_2)t^2|\Delta u|^2= 5C_2 t^2|\Delta u|^2$, and 
\begin{align}
C_5&= C(n, K, \mathscr{M})(1+ T) \Big[1+ T(1+ \frac{1}{\rho^2})\Big] \label{2.22} \\
C_6&= C(n, K, \mathscr{M}) \Big[1+ T(1+ \frac{1}{\rho^2})\Big]^{-2} \label{2.23}
\end{align}
Define $\gamma(x)= \phi(\frac{d(x, p)}{\rho})$, then $\gamma(x)F(x, t)$ attains its maximum at a point $(x_0, t_0)\in B_{p}(2\rho)\times [0, T]$. The rest computation is at $(x_0, t_0)$, 
\begin{align}
0\leq \Big(\frac{\partial}{\partial t}- \Delta \Big) (\gamma F)\leq \gamma \Big(-\frac{C_6}{t}F^2+ \frac{C_5}{t}\Big)- \Delta \gamma \cdot F- 2\nabla \gamma \nabla F \nonumber
\end{align}
Note at $(x_0, t_0)$, $\nabla(\gamma F)= 0$, let $G= (\gamma F)(x_0, t_0)$, we get
\begin{align}
0\leq -\frac{C_6}{t}G^2+ \Big(2\frac{|\nabla \gamma|^2}{\gamma}- \Delta \gamma\Big)G+ \frac{C_5}{t} \label{2.24}
\end{align}
and 
\begin{align}
\Big(2\frac{|\nabla \gamma|^2}{\gamma}- \Delta \gamma\Big)&= \frac{2}{\rho^2}\cdot \frac{|\phi '|^2}{\phi}- \frac{\phi ''}{\rho^2}- \frac{\phi '}{\rho} \Delta d(x, p) \nonumber \\
&\leq \frac{32}{\rho^2}+ \frac{2}{\rho}\cdot \coth \Big(\sqrt{\frac{K}{n- 1}}\rho\Big) \label{2.25}
\end{align}
in the last inequality we used (\ref{2.19}), $Rc\geq -Kg$ and Laplacian comparison theorem. From (\ref{2.22})-(\ref{2.25}), 
\begin{align}
0\leq -G^2&+ C(n, K, \mathscr{M})\Big[1+ T(1+ \frac{1}{\rho^2})\Big]^2T\cdot \Big[\frac{1}{\rho^2}+ \frac{1}{\rho}\coth \Big(\sqrt{\frac{K}{n- 1}}\rho\Big)\Big] G \nonumber \\
&+ C(n, K, \mathscr{M}) \Big[1+ T(1+ \frac{1}{\rho^2})\Big]^3 (1+ T) \nonumber
\end{align}
then it is easy to get
\begin{align}
G&\leq C(n, K, \mathscr{M})\cdot\Big[1+ T(1+ \frac{1}{\rho^2})\Big]^2(1+ T) \nonumber \\
&\quad \cdot \Big[\Big(\frac{1}{\rho^2}+ \frac{1}{\rho}\Big)\coth \Big(\sqrt{\frac{K}{n- 1}}\rho\Big) + 1+ T(1+ \frac{1}{\rho^2})\Big]  \nonumber
\end{align}
Hence on $B_{p}(\rho)$,
\begin{align}
t^2|\Delta u|^2&\leq C_3^{-1}F\leq C_3^{-1} G \nonumber \\
&\leq C(n, K, \mathscr{M})\cdot\Big[1+ T(1+ \frac{1}{\rho^2})\Big]^2\cdot \Big[\Big(\frac{1}{\rho^2}+ 1\Big)\cdot \Big(T+ \coth \Big(\sqrt{\frac{K}{n- 1}}\rho\Big) \Big)+ 1\Big]  \nonumber
\end{align}
Taking square root in the above inequality, we can get our conclusion. 
}
\qed

Let $\rho\rightarrow \infty$, we get the following global estimate.

\begin{cor}\label{cor 2.3.1}
{Suppose $(M^n, g)$ is a complete Riemannian manifold with $Rc\geq -Kg$, and $|u(x, t)|\leq \mathscr{M}$ is a solution to the heat equation on $M^n\times [0, T]$, where $K$, $\mathscr{M}$, $T$ are positive constants. Then 
\begin{equation}\label{2.3.1.1}
{t|\Delta u|\leq  C(n, K, \mathscr{M})(1+ T)^2
}
\end{equation}
holds on $M^n\times [0, T]$.
}
\end{cor}

We also need a maximum principle due originally to Karp and Li (see \cite{KL}) whose statement can be found (in more generalized form) in Ni and Tam's paper \cite{NT} (see Theorem 1.2 there). The statement of this maximum principle is as follows.

\begin{theorem}[Karp-Li and Ni-Tam]\label{thm 2.3.1}
{Suppose $(M^n, g)$ is a complete Riemannian manifold and $h(x, t)$ is a smooth function on $M^n\times [0, T]$ such that $\Big(\frac{\partial}{\partial t}- \Delta\Big) f(x, t)\leq 0$ whenever $f(x, t)\geq 0$. Assume that 
\begin{equation}\nonumber
{\int_{0}^T\int_{M^n} e^{-a\cdot d^2(x, p)}f_{+}^2(x, s)d\mu(x) ds< \infty
}
\end{equation}
for some $a> 0$, where $p$ is a fixed point on $M^n$ and $f_{+}(x, t)\doteqdot \max\{f(x, t), 0\}$. If $f(x, 0)\leq 0$ for all $x\in M^n$, then $f(x ,t)\leq 0$ for all $(x, t)\in M^n\times [0, T]$.
}
\end{theorem}

Now we are ready to prove Hamilton's Theorem in the complete case.
\begin{theorem}\label{thm 2.4}
{Suppose $(M^n, g)$ is a complete Riemannian manifold with $Rc\geq -Kg$, and $0< u(x, t)\leq \mathscr{M}$ is a solution to the heat equation on $M^n\times [0, T]$, where $K$, $\mathscr{M}$, $T$ are positive constants. Then 
\begin{equation}\label{2.4.1}
{t\Big(\frac{\Delta u}{u}+ \frac{|\nabla u|^2}{u^2}\Big)\leq n+ (4+ 2Kt)\ln\Big(\frac{\mathscr{M}}{u}\Big)
}
\end{equation}
}
\end{theorem}

\pf
{Define $u_{\epsilon}= u+ \epsilon$ for $\epsilon> 0$, we obtain a solution satisfying $\epsilon< u_{\epsilon}\leq \mathscr{M}+ \epsilon\doteqdot \mathscr{M}_{\epsilon}$. Once the estimate has been proved for $u_{\epsilon}$, the theorem will follow by letting $\epsilon \rightarrow 0$. Consider the function
\begin{equation}\nonumber
{F(x, t)= t\Big(\Delta u_{\epsilon}+ \frac{|\nabla u_{\epsilon}|^2}{u_{\epsilon}}\Big)- u_{\epsilon}\Big[ n+ (4+ 2Kt)\ln\Big(\frac{\mathscr{M}_{\epsilon}}{u_{\epsilon}}\Big) \Big]
}
\end{equation}
A long but straightforward computation gives
\begin{align}
\Big(\frac{\partial}{\partial t}- \Delta\Big) F(x, t)&= u_{\epsilon} \Big[- 2t|\nabla^2 \ln u_{\epsilon}|^2+ \Delta \ln u_{\epsilon}- (2+ 2Kt) |\nabla \ln u_{\epsilon}|^2   \nonumber \\
&\quad \quad - 2tRc(\nabla \ln u_{\epsilon}, \nabla \ln u_{\epsilon})- 2K\ln\Big(\frac{\mathscr{M}_{\epsilon}}{u_{\epsilon}}\Big)\Big] \nonumber \\
&\leq u_{\epsilon} \Big[-\frac{2t}{n}|\Delta \ln u_{\epsilon}|^2+ \Delta \ln u_{\epsilon}- 2|\nabla \ln u_{\epsilon}|^2 \Big] \label{2.4.2}
\end{align}
If $F(x, t)\geq 0$ at $(x, t)$, then 
{\begin{equation}\label{2.4.3}
{-2|\nabla \ln u_{\epsilon}|^2\leq \Delta \ln u_{\epsilon}- \frac{n}{t}
}
\end{equation}
From (\ref{2.4.2}) and (\ref{2.4.3}), 
\begin{align}
\Big(\frac{\partial}{\partial t}- \Delta\Big) F(x, t)&=\leq  u_{\epsilon} \Big[-\frac{2t}{n}|\Delta \ln u_{\epsilon}|^2+ 2\Delta \ln u_{\epsilon}- \frac{n}{t} \Big] \nonumber \\
&\leq -\frac{n}{2t}< 0 \label{2.4.4}
\end{align}

In (\ref{2.21}) let $\rho\rightarrow \infty$, 
\begin{equation}\label{2.4.5}
{t|\nabla u|^2\leq C(K, \mathscr{M}_{\epsilon}, T)
}
\end{equation}

From (\ref{2.4.5}) and (\ref{2.3.1.1}), 
\begin{align}
F_{+}^2(x, t)\leq \Big[t\Big(\Delta u_{\epsilon}+ \frac{|\nabla u_{\epsilon}|^2}{u_{\epsilon}}\Big)\Big]^2\leq C(\epsilon, n, K, \mathscr{M}_{\epsilon}, T) \label{2.4.6}
\end{align}

Using (\ref{2.4.6}), for any $p\in M^n$ and $\rho> 0$,
\begin{align}
&\int_{0}^{T}\int_{B_p(\rho)} \exp{\Big(-d^2(x, p)\Big)}F_{+}^2(x, t) d\mu(x) dt \nonumber \\
&\leq C(\epsilon, n, K, \mathscr{M}_{\epsilon}, T) \int_{M^n} \exp{\Big(-d^2(x, p)\Big)} d\mu(x) \leq C \label{2.4.7}
\end{align}
in the last equality we used volume comparison theorem and $Rc\geq -Kg$. Let $\rho\rightarrow \infty$, 
\begin{equation}\label{2.4.8}
{\int_{0}^{T}\int_{M^n} \exp{-d^2(x, p)}F_{+}^2(x, t) d\mu(x) dt \leq C< \infty
}
\end{equation}

From (\ref{2.4.4}) and (\ref{2.4.8}), using Theorem \ref{thm 2.3.1}, we get $F(x, t)\leq 0$ for all $0\leq t\leq T$. The conclusion is proved.
}
\qed

The following corollary gives the upper bound of Laplacian of the heat kernel.
\begin{cor}\label{cor 2.5}
{Suppose $(M^n, g)$ is a complete Riemannian manifold with $Rc\geq -Kg$, and $H(x, y, t)$ is the heat kernel on $M^n$, and $0< t\leq T$, where $K$, $T$ are positive constants. Then
\begin{align}
&\Big(\Delta H+ \frac{|\nabla H|^2}{H}\Big)(x, y, t) \nonumber \\
&\leq \frac{2H(x, y,t)}{t}\Big\{n+ (4+ Kt)\ln \Big[\frac{C(K, T)}{H(x, y, t) V^{\frac{1}{2}}(B_{x}(\sqrt{\frac{t}{2}})) V^{\frac{1}{2}}(B_{y}(\sqrt{\frac{t}{2}}))}\Big]\Big\} \label{2.5.1}
\end{align}
}
\end{cor}

\pf
{Note if $s\in [\frac{t}{2}, t]$, from \cite{LY},
\begin{align}
H(x, y, t)\leq C(K, T)\cdot V^{\frac{1}{2}}\Big(B_{x}\Big(\sqrt{\frac{t}{2}}\Big) \Big) V^{\frac{1}{2}}\Big( B_{y}\Big(\sqrt{\frac{t}{2}}\Big) \Big) \nonumber 
\end{align}
Then apply Theorem \ref{thm 2.4} on $u(x, s)= H(x, y, s+ \frac{t}{2})$ and $M^n\times[0, \frac{t}{2}]$, conclusion follows from (\ref{2.4.1}).
}
\qed


\begin{thebibliography}{99}
\bibitem{CLY} Siu Yuen, Cheng; Peter, Li; Shing Tung, Yau:
\emph{On the upper estimate of the heat kernel of a complete Riemannian manifold}. Amer. J. Math. 103 (1981), no. 5, 1021-1063. 

\bibitem{RFAA} Chow, Bennett; Chu, Sun-Chin; Glickenstein, David; Guenther, Christine; Isenberg, James; Ivey, Tom; Knopf, Dan; Lu, Peng; Luo, Feng; Ni, Lei:
\emph{The Ricci flow: techniques and applications. Part II. Analytic aspects}. Mathematical Surveys and Monographs, 144. American Mathematical Society, Providence, RI, 2008. 458 pp.

\bibitem{RFTA} Chow, Bennett; Chu, Sun-Chin; Glickenstein, David; Guenther, Christine; Isenberg, James; Ivey, Tom; Knopf, Dan; Lu, Peng; Luo, Feng; Ni, Lei:
\emph{The Ricci flow: techniques and applications. Part III. Geometric-analytic aspects}. Mathematical Surveys and Monographs, 163. American Mathematical Society, Providence, RI, 2010. 517 pp.

\bibitem{GL} Nicola, Garofalo \& Ermanno, Lanconelli:
\emph{Asymptotic behavior of fundamental solutions and potential theory of parabolic operators with variable coefficients}. Math. Ann. 283 (1989), no. 2, 211-239.

\bibitem{Gri} Alexander, Grigor'yan:
\emph{Heat kernel and analysis on manifolds}. AMS/IP Studies in Advanced Mathematics, 47. American Mathematical Society, Providence, RI; International Press, Boston, MA, 2009. xviii+482 pp.

\bibitem{Ham} Richard S. Hamilton:
\emph{A matrix Harnack estimate for the heat equation}. Comm. Anal. Geom. 1 (1993), no. 1, 113-126.

\bibitem{KL} Leon, Karp \& Peter, Li:
\emph{The heat equation on complete riemannian manifolds}. Unpublished paper, 1982.

\bibitem{Kots} Brett L. Kotschwar:
\emph{Hamilton's gradient estimate for the heat kernel on complete manifolds}. Proc. Amer. Math. Soc. 135 (2007), no. 9, 3013-3019 

\bibitem{Li}Peter, Li:
\emph{Geometric Analysis}. Cambridge Studies in Advanced Mathematics, 134. Cambridge University Press, 2012, 416 pp.

\bibitem{LY} Peter, Li \& Shing-Tung, Yau:  
\emph{On the parabolic kernel of the Schrödinger operator}. Acta Math. 156 (1986), no. 3-4, 153-201. 

\bibitem{MS} Paul Malliavin \& Daniel W. Stroock:
\emph{Short time behavior of the heat kernel and its logarithmic derivatives}, J. Differential Geom. 44(3)
(1996), 550-570.

\bibitem{Ni} Lei, Ni:
\emph{The entropy formula for linear heat equation}. J. Geom. Anal. 14 (2004), no. 1, 87-100
 
\bibitem{Niadd} Lei, Ni:
\emph{Addenda to: "The entropy formula for linear heat equation''}. J. Geom. Anal. 14 (2004), no. 1, 87-100; MR2030576]. J. Geom. Anal. 14 (2004), no. 2, 369-374. 

\bibitem{LYH} Lei, Ni:
\emph{A note on Perelman's LYH-type inequality}. Comm. Anal. Geom. 14 (2006), no. 5, 883-905.

\bibitem{Nilarge} Lei, Ni:
\emph{The large time asymptotics of the entropy}. Complex analysis, 301-306, Trends Math., Birkhäuser/Springer Basel AG, Basel, 2010.
 
\bibitem{NT} Lei, Ni \& Luen-Fai, Tam:
\emph{K\"ahler-Ricci flow and the Poincaré-Lelong equation}. Comm. Anal. Geom. 12 (2004), no. 1-2, 111-141.
 
\bibitem{Pere} G. Perelman:
\emph{The entropy formula for the Ricci flow and its geometric applications}. 2002, arXiv:math.DG/0211159.

\bibitem{Shi} Wan-Xiong, Shi:
\emph{Deforming the metric on complete Riemannian manifolds}. J. Differential Geom. 30 (1989), no. 1, 223-301. 


\bibitem{Va} S. R. S. Varadhan:
\emph{On the behavior of the fundamental solution of the heat equation with variable coefficients}. 
Comm. Pure Appl. Math. 20, 1967, 431-455.

\bibitem{Wu} Jia-Yong, Wu:
\emph{Sharp Hamilton's Laplacian estimate for the heat kernel on complete manifolds}, to appear on Proceedings of AMS.

\end{thebibliography}
\end{document}